\documentclass[review]{siamart1116}

\usepackage{latexsym}
\usepackage{amsmath}
\usepackage{amssymb}
\usepackage{amsfonts}
\usepackage{verbatim}
\usepackage{mathrsfs}
\usepackage{amsfonts}
\usepackage{graphicx}
\usepackage{colortbl,dcolumn}
\usepackage{amsmath}
\usepackage{psfrag}
\usepackage{booktabs}

\usepackage{mathtools,bbm}

\newtheorem{tm}{Theorem}[section]

\newtheorem{ap}{Assumption}[section]

\newtheorem{prop}{Proposition}[section]
\newtheorem{lm}{Lemma}[section]
\newtheorem{cor}{Corollary}[section]

\newcommand{\N}{\mathbb N}

\newcommand{\n}{\mathbf n}
\newcommand{\E}{\mathbf E}
\renewcommand{\H}{\mathbf H}
\newcommand{\W}{\mathbf W}

\renewcommand{\S}{\mathbf S}

\newcommand{\R}{\mathbb R}
\renewcommand{\P}{\mathbb P}
\newcommand{\EE}{\mathbb E}
\newcommand{\FF}{\mathbb F}
\newcommand{\GG}{\mathbb G}
\newcommand{\UU}{\mathbb U}
\newcommand{\VV}{\mathbb V}

\renewcommand{\AA}{\mathcal A}

\newcommand{\LL}{\mathcal L}
\newcommand{\OO}{\mathcal O}
\newcommand{\HH}{\mathcal H}

\newcommand{\FFF}{\mathscr F}
\newcommand{\OOO} {\mathscr O}
\newcommand{\<}{\langle}
\renewcommand{\>}{\rangle}
\newcounter{RomanNumber}
\newcommand{\MyRoman}[1]{\setcounter{RomanNumber}{#1}
\rm\Roman{RomanNumber}}
\newcommand{\Tr}{\mathrm{Tr}}
\newcommand{\ddiv}{\mathrm{div}}
\newcommand{\ccurl}{\mathrm{curl}}

\newcommand{\TheTitle}{Exponential Integrators for Stochastic Maxwell's Equations Driven by It\^o Noise} 
\newcommand{\TheAuthors}{David Cohen, Jianbo Cui, Jialin Hong and Liying Sun}

\headers{}{\TheAuthors}

\title{{\TheTitle}\thanks{Submitted to the editors in DATE.
\funding{his work was supported by the 
National Natural Science Foundation of China 
(NO. 91530118, NO. 91130003, NO. 11021101, NO. 91630312 and NO. 11290142), 
the Swedish Foundation for International Cooperation in Research and Higher Education 
(STINT project nr. $CH2016-6729$), as well as the Swedish Research Council (VR) 
(projects nr. $2013−4562$ and $2018-04443$). 
The computations were performed on resources provided by the Swedish National Infrastructure 
for Computing (SNIC) at HPC2N, Ume{\aa} University.}}}

\author{
David Cohen
\thanks{Department of Mathematics and Mathematical Statistics, Ume$\mathring{a}$ University, 90187 Ume$\mathring{a}$, Sweden
(\email{david.cohen@umu.se})}
\and
Jianbo Cui 
\thanks{1. LSEC, ICMSEC, 
Academy of Mathematics and Systems Science, Chinese Academy of Sciences, Beijing,  100190, China\qquad
2. School of Mathematical Science, University of Chinese Academy of Sciences, Beijing, 100049, China 
(\email{jianbocui@lsec.cc.ac.cn}, \email{hjl@lsec.cc.ac.cn},\email{liyingsun@lsec.cc.ac.cn}(corresponding author))}
\and 
Jialin Hong \footnotemark[3]
\and 
Liying Sun \footnotemark[3]
}

\usepackage{amsopn}

\ifpdf
\hypersetup{
  pdftitle={\TheTitle},
  pdfauthor={\TheAuthors}
}
\fi

\usepackage{ulem}
\usepackage{subfigure}
\usepackage{graphicx}
\usepackage{enumerate}
\usepackage{amsmath,bm}
\usepackage{latexsym}
\usepackage{float}

\begin{document}

\maketitle

\begin{abstract}
This article  presents explicit exponential integrators for stochastic Maxwell's equations driven by both multiplicative and additive noises. 
By utilizing the regularity estimate of the mild solution,
we first prove that the strong order of the numerical approximation is $\frac 12$ for general multiplicative noise. Combing a proper decomposition with the stochastic Fubini's theorem,
the strong order of the proposed scheme  is shown to be $1$ for additive noise. 
Moreover, for linear stochastic Maxwell's equation 
with additive noise, the proposed time integrator is shown to preserve exactly 
the symplectic structure, the evolution of the energy as well as the evolution 
of the divergence in the sense of expectation. Several numerical experiments 
are presented in order to verify our theoretical findings.
\end{abstract}


\begin{AMS}
{60H35}, {60H15, 35Q61.}
\end{AMS}

\section{Introduction}
In the context of electromagnetism, a common way to model precise 
microscopic origins of randomness (such as thermal 
motion of electrically charged micro-particles) 
is by means of stochastic Maxwell's equations \cite{RKT89}. 
Further applications of stochastic Maxwell's equations are: 
In \cite{OGN96}, a stochastic model of Maxwell's field equations in $1+1$ dimension 
is shown to be a simple modification of a random walk model due to Kac, 
which provides a basis for the telegraph equations. 
The work \cite{KS14} studies the propagation of ultra-short solitons 
in a cubic nonlinear medium modeled by nonlinear Maxwell's equations 
with stochastic variations of media. 
To simulate a coplanar waveguide with uncertain material parameters, 
time-harmonic Maxwell's equations are considered in \cite{BS15}. 
For linear stochastic Maxwell's equations driven by additive noise,  
the work \cite{HJZ14} proves that the problem is a stochastic 
Hamiltonian  partial differential equation whose phase flow preserves 
the multi-symplectic geometric structure. 
In addition, the averaged energy along the flow increases linearly with respect to time and 
the flow preserves the divergence in the sense of expectation, see \cite{CHZ16}.
Let us finally mention that linear stochastic Maxwell's equations are relevant 
in various physical applications, see e.g. \cite[Chapter 3]{RKT89}. 

We now review the literature on the numerical discretisation of stochastic Maxwell's equations. 
The work \cite{Zha08} performs a numerical analysis of the finite element method and 
discontinuous Galerkin method for stochastic Maxwell's equations driven by colored noise. 
A stochastic multi-symplectic method for $3$ dimensional problems with additive noise, 
based on stochastic variational principle, is studied in \cite{HJZ14}. 
In particular, it is shown that the implicit numerical scheme preserves 
a discrete stochastic multi-symplectic conservation law. 
The work \cite{CHZ16} inspects geometric properties of the stochastic Maxwell's equation with additive noise, 
namely the behavior of averaged energy and divergence, see below for further details. 
Especially, the authors of \cite{CHZ16} investigate three novel stochastic multi-symplectic 
(implicit in time) methods preserving discrete versions of the averaged divergence. 
None of the proposed numerical schemes exactly preserve the behavior of the averaged energy. 
The work \cite{HJZC17} proposes a stochastic multi-symplectic wavelet collocation method 
for the approximation of stochastic Maxwell's equations with multiplicative noise (in the Stratonovich sense). 
For the same stochastic Maxwell's equation as the one considered in this paper 
(see below for a precise definition), 
the recent reference \cite{CHJ18} shows that the backward Euler--Maruyama method 
converges with mean-square convergence rate $\frac 12$. Finally, the preprint \cite{CHJ18a} 
studies implicit Runge--Kutta schemes for stochastic Maxwell's equation with additive noise. 
In particular, a mean-square convergence of order $1$ is obtained. 

In the present paper, we construct and analyse an exponential integrator 
for stochastic Maxwell's equations which is explicit (thus computationally more efficient 
than the above mentioned time integrators) and which enjoys excellent long-time behavior.  
Observe that exponential integrators are widely used for efficient 
time integrations of deterministic differential equations,  
see for instance \cite{hls98,MR2413146,ho10,cg12} 
and more specially \cite{tk02,ntb07,MR2437586,MR2528933,paz13} and references therein for Maxwell-type equations. 
In recent years, exponential integrators have been analysed 
in the context of stochastic (partial) differential equations (S(P)DEs). 
Without being too exhaustive, we mention analysis and applications of 
such numerical schemes for the following problems: 
stochastic differential equations \cite{sxz12,kb14,kcb17}; 
stochastic parabolic equations \cite{jk09,lt13,BCH18,CH18,acqs18}; 
stochastic Schr\"odinger equations \cite{CA17,CD17,CHLZ16}; 
stochastic wave equations \cite{cls13,W15,cq16,ACLW16,QW17}
and references therein.

The main contributions of the present paper are: 
\begin{itemize}
\item a strong convergence analysis of an explicit exponential 
integrator for stochastic Maxwell's equations in $\mathbb R^3$. 
By making use of regularity estimates of the exact and numerical solutions, 
the strong convergence order is shown to be $\frac 12$ for general multiplicative noise. 
Furthermore, by using a proper decomposition and stochastic Fubini's theorem, 
we prove that the strong convergence order of the proposed scheme can achieve $1$.
\item an analysis of long-time conservation properties of an explicit exponential integrator  
for linear stochastic Maxwell's equations driven by additive noise. 
Especially, we show that the proposed explicit time integrator is symplectic and 
satisfies a trace formula for the energy for all times, 
i.\,e. the linear drift of the averaged energy is preserved for all times. 
In addition, the numerical solution preserves 
the averaged divergence. This shows that the exponential integrator inherits 
the geometric structure and the dynamical behavior  of the flow of the 
linear stochastic Maxwell's equations. This is not the case for classical time integrators 
such as Euler--Maruyama type schemes. 
\item an efficient numerical implementation of two-dimensional models of stochastic Maxwell’s equations 
by explicit time integrators. 
\end{itemize}
We would like to remark that the proofs of strong convergence 
for the exponential integrator use similar ideas present 
in various proofs of strong convergence from the literature. 
But, to the best of our knowledge, the present paper offers the first explicit time integrator 
for linear stochastic Maxwell's equations that is of strong order $1$, symplectic, exactly preserves 
the linear drift of the averaged energy, and preserves the averaged divergence for all times. 
A weak convergence analysis of the proposed scheme for stochastic Maxwell's equations 
driven by multiplicative noise will be reported elsewhere. 

An outline of the paper is as follows. 
Section~\ref{sec;2} sets notations and introduces the stochastic Maxwell's equation.  
This section also presents assumptions to guarantee existence and uniqueness 
of the exact solution to the problem and shows its H\"older continuity. 
The exponential integrator for stochastic Maxwell's equation is introduced in Section~\ref{sect-EXP}, 
where we also prove its strong order of convergence for additive and multiplicative noise. 
In Section~\ref{sect-LSM}, we show that the proposed scheme has several interesting geometric properties: 
it preserves the evolution laws of the averaged energy, the evolution laws of the divergence, 
and the symplectic structure of the original linear stochastic Maxwell's equations 
with additive noise. We conclude the paper by presenting numerical experiments 
supporting our theoretical results in Section~\ref{sect-NE}.


\section{Well-posedness of stochastic Maxwell's equations}\label{sec;2}
We consider the stochastic Maxwell's equation driven by multiplicative It\^o noise 
\begin{align}
\label{mod;CMAX}
\begin{split}
\text d\UU&=A\UU\, \text dt+\FF(\UU)\,\text dt+\GG(\UU)\,\text dW, \quad {t} \;{\in}\; (0,+\infty),
\\
\UU(0)&=(\E_0^\top,\H_0^\top)^\top
\end{split}
\end{align}
supplemented with the boundary condition of a perfect conductor $\n\times\E=0$ as in \cite{HJZ14}.
Here, $\UU=(\E^\top,\H^\top)^\top$, is $\R^6$-valued function whose domain $\OO$ is a bounded and simply connected domain in $\R^3$ 
with smooth boundary $\partial \OO$. 
The unit outward normal vector to $\partial \OO$ is denoted by $\n$. 
Moreover, $\text dW$ stands for the formal time derivative of a $Q$-Wiener process $W$ on a stochastic 
basis $(\Omega, \FFF, \{\FFF_t\}_{t\ge 0}, \P)$. The $Q$-Wiener process can be written as 
$W({\bf x},t)=\sum\limits_{k\in \N_+} Q^{\frac12}e_k({\bf x})\beta_k(t)$, where 
$\{\beta_k\}_{k\in \N_+}$ is a sequence of mutually independent and identically distributed 
$\R$-valued standard Brownian motions; $\{e_k\}_{k\in \N_+}$ is an orthonormal basis of $U:=\LL^2(\OO; \R)$ 
consisting of eigenfunctions of a symmetric, nonnegative and of finite trace linear operator $Q$, i.\,e., 
$Q e_k=\eta_k e_k$, with $\eta_k\ge 0$ for $k\in \N_+$. 
Assumptions on $\FF$ and $\GG$ are provided below.

The Maxwell's operator $A$ is defined by 
\begin{align}\label{Mop}
A\begin{pmatrix}
\E\\ \H
\end{pmatrix}:=
\begin{pmatrix}
0 & \epsilon^{-1}\nabla\times \\ -\mu^{-1}\nabla\times & 0
\end{pmatrix}
\begin{pmatrix}
\E\\ \H
\end{pmatrix}
=
\begin{pmatrix}
\epsilon^{-1}\nabla\times\H\\ -\mu^{-1}\nabla\times\E
\end{pmatrix}.
\end{align}
It has the domain $D(A):=H_0(\ccurl,\OOO)\times H(\ccurl, \OOO)$, where
\begin{align*}
&H(\ccurl,\OOO):=\{\mathbf U\in (\LL^2(\OOO))^3:\nabla\times \mathbf U\in (\LL^2(\OOO))^3\},
\end{align*}
is termed by the {${\ccurl}$}-space and
\begin{align*}
&H_0(\ccurl,\OOO):=\{\mathbf U\in H(\ccurl,\OOO):\n\times \mathbf U|_{\partial \OOO}=\mathbf 0\}
\end{align*}
is the subspace of $H(\ccurl,\OOO)$ with zero tangential trace.
In addition, $\epsilon$ and $\mu$ are bounded and uniformly positive definite functions: 
$$
\epsilon,\mu\in \LL^\infty(\OOO),\quad \epsilon,
\mu\geq \kappa>0
$$
with $\kappa$ being a positive constant. These conditions on $\epsilon,\mu$ ensure 
that the Hilbert space $V:=(\LL^2(\OOO))^3\times (\LL^2(\OOO))^3$ is equipped 
with the weighted scalar product
\begin{align*}
\left\<
\begin{pmatrix}
\E_1 \\ \H_1
\end{pmatrix},
\begin{pmatrix}
\E_2 \\ \H_2
\end{pmatrix}
\right\>_V
=\int_\OOO\left( \mu\<\H_1,\H_2\>+\epsilon\<\E_1,\E_2\>\right)\,\text d{\bf x},
\end{align*}
where $\<\cdot,\cdot\>$ stands for the standard Euclidean inner product.
This weighted scalar product is equivalent to the standard inner product on $(\LL^2(\OO))^6$. 
Moreover, the corresponding norm, which stands for the electromagnetic energy 
of the physical system, induced by this inner product reads 
\begin{align*}
\left\|
\begin{pmatrix}
\E \\ \H
\end{pmatrix}
\right\|_V^2
=\int_\OOO \left(\mu\|\H\|^2+\epsilon\|\E\|^2\right)\,\text d{\bf x}
\end{align*}
with $\|\cdot\|$ being the Euclidean norm.
Based on the norm $\|\cdot\|_V$, the associated graph norm of $A$ is defined by 
\begin{align*}
\|\VV\|_{D(A)}^2:=\|\VV\|_V^2+\|A\VV\|_V^2.
\end{align*} 
It is well known that Maxwell's operator $A$ is closed and that $D(A)$ equipped
with the graph norm is a Banach space, see e.g. \cite{Mon03}. 
Moreover, $A$ is skew-adjoint, in particular, for all $\VV_1,\VV_2\in D(A)$, 
\begin{align*}
\<A\VV_1,\VV_2\>_V=-\<\VV_1,A\VV_2\>_V.
\end{align*} 
In addition, the operator $A$ generates a unitary $C_0$-group $\S(t):=\exp (tA)$ via Stone's theorem, 
see for example \cite{HJS15}. 
According to the definition of unitary groups, one has
\begin{align}
\label{unigro}
\|\S(t)\VV\|_V=\|\VV\|_V\quad\text{for all}\quad\VV\in V,
\end{align}
which means that the electromagnetic energy is preserved, for Maxwell's operator, see \cite{HP15}.
Besides, the unitary group $\S(t)$ satisfies the following properties which will be made use of in the next section.
\begin{lm}[Theorem~3 with $\mathbf{q=0}$ in \cite{BT79}]
\label{lm;SG}
For the semigroup $\{\S(t);t\geq 0\}$ on $V$, it holds that
\begin{align}
\|\S(t)-Id\|_{L\left(D(A);V\right)}\leq Ct, 
\end{align}
where the constant $C$ does not depend on $t$. 
Here, $L(D(A);V)$ denotes the space of bounded linear operators from $D(A)$ to $V$.
\end{lm}
Observe that, throughout the paper, $C$ stands for a constant that may vary from line to line. 

For two real-valued separable Hilbert spaces $(H_1,\<\cdot,\cdot\>_{H_1},\|\cdot\|_{H_1})$ 
and $(H_2,\<\cdot,\cdot\>_{H_2},\|\cdot\|_{H_2})$, we denote the set 
of Hilbert--Schmidt operators from $H_1$ to $H_2$ by $\LL_2(H_1,H_2)$. It will be 
equipped with the norm
$$
\|\Gamma\|^2_{\LL_2(H_1,H_2)}:=\sum\limits_{i=1}^\infty\|\Gamma\phi_i\|_{H_2}^2, 
$$
where $\{\phi_i\}_{i\in\N_+}$ is any orthonormal basis of $H_1$. 
Furthermore, let $Q^\frac 12$ be the unique positive square root of the linear operator $Q$ 
(defining the noise $W$). 
We also introduce the separable Hilbert space $U_0:=Q^{\frac 12}U$ endowed with the inner product 
$\<u_1,u_2\>_{U_0}:=\<Q^{-\frac 12}u_1,Q^{-\frac 12}u_2\>_U$ for $u_1,u_2\in U_0$, 
where we recall that $U=\LL^2(\OO; \R)$. 

\begin{lm}
As a consequence of Lemma~\ref{lm;SG}, for any $\Phi\in \LL_2\left(U_0,D(A)\right)$ and any $t\geq 0,$ we have
\begin{align}
\label{SGD}
\|\left(\S(t)-Id\right)\Phi\|_{\LL_2(U_0,V)}\leq Ct\|\Phi\|_{\LL_2(U_0,D(A))}.
\end{align}
\end{lm}
{\bf Proof} Thanks to Lemma~\ref{lm;SG} and the definition of the Hilbert--Schmidt norm, we know that, 
for $\{e_k\}_{k\in\N_+}$  an orthonormal basis of $U$, 
\begin{align*}
\|\left(\S(t)-Id\right)\Phi\|_{\LL_2(U_0,V)}^2
&=\sum\limits_{k\in \N_+} \|\left(\S(t)-Id\right)\Phi Q^\frac 12e_k\|^2_V\\
&\leq Ct^2\sum\limits_{k\in \N_+} \|\Phi Q^\frac 12e_k\|^2_{D(A)}\leq Ct^2\|\Phi\|_{\LL_2(U_0,D(A))}^2,
\end{align*}
which proves the claim.
\hfill$\square$

To guarantee existence and uniqueness of strong solutions to \eqref{mod;CMAX}, we make the following assumptions: 

\begin{ap}[Coefficients]
\label{ap;1}
Assume that the coefficients of Maxwell's operator \eqref{Mop} satisfy
$$
\epsilon,\mu\in \LL^\infty(\OOO),\quad \epsilon,
\mu\geq \kappa>0
$$
with some positive constant $\kappa$.
\end{ap}

\begin{ap}[Initial value]
\label{ap;4}
The initial value $\UU(0)$ of the stochastic Maxwell's equation \eqref{mod;CMAX} is a $D(A)$-valued stochastic process with 
$\EE\left[\|\UU(0)\|^p_{D(A)}\right]<\infty$ for any $p\ge 1$.
\end{ap}

\begin{ap}[Nonlinearity]
\label{ap;2}
We assume that the operator $\FF\colon V\to V$ is
continuous and that there exists constants $C_\FF,C_\FF^1>0$ such that 
\begin{align*}
&\|\FF(\VV_1)-\FF(\VV_2)\|_V\leq C_\FF\|\VV_1-\VV_2\|_V, \quad \VV_1,\VV_2\in V,\\
&\|\FF(\VV_1)-\FF(\VV_2)\|_{D(A)}\leq C_\FF^1\|\VV_1-\VV_2\|_{D(A)}, \quad \VV_1,\VV_2\in D(A),\\
&\|\FF(\VV)\|_V\leq C_\FF(1+\|\VV\|_V),\quad \VV\in V,\\
&\|\FF(\VV)\|_{D(A)}\leq C_\FF^1\left(1+\|\VV\|_{D(A)}\right),\quad \VV\in D(A).
\end{align*}
\end{ap}

\begin{ap}[Noise]
\label{ap;3}
We assume that the operator $\GG\colon V\to \LL_2(U_0,V)$ satisfies 
\begin{equation}
\begin{split}
\label{con;G}
&\|\GG(\VV_1)-\GG(\VV_2)\|_{\LL_2(U_0,V)}\leq C_\GG\|\VV_1-\VV_2\|_V,
\; \VV_1,\VV_2\in V,\\
&\|\GG(\VV_1)-\GG(\VV_2)\|_{\LL_2(U_0,{D(A)})}\leq C_\GG^1\|\VV_1-\VV_2\|_{D(A)},\;  \VV_1, \VV_2\in D(A),\\
&\|\GG(\VV)\|_{\LL_2(U_0,V)}\leq C_\GG(1+\|\VV\|_V),\quad \VV\in V,\\
&\|\GG(\VV)\|_{\LL_2(U_0,D(A))}\leq C_\GG^1(1+\|\VV\|_{D(A)}),\quad \VV\in D(A),
\end{split}
\end{equation}
where $C_\GG,C_\GG^1>0$  may depend on the operator $Q$. 
We recall that $\LL_2(U_0,V)$ and $\LL_2(U_0,D(A))$ denote 
the spaces of Hilbert--Schmidt operators from $U_0$ to $V$, resp. to $D(A)$.
\end{ap}

We now present two examples of an operator $\GG$ verifying Assumption~\ref{ap;3} 
(we only prove one of the inequality in \eqref{con;G}, the others follow in a similar way).

For the first example (inspired by \cite{HJZ14}), let $\OO=[0,1]^3$, $\epsilon=\mu=1$ 
and consider $\GG\equiv(\lambda_1 ,\lambda_1,\lambda_1,\lambda_2 ,\lambda_2 ,\lambda_2)^T$ 
for two real numbers $\lambda_1$ and $\lambda_2$. 
The stochastic Maxwell's equation \eqref{mod;CMAX} then becomes an SPDE 
driven by additive noise. In this case, one chooses the orthonormal basis of $U$ 
to be $\sin(i\pi x_1)\sin(j\pi x_2)\sin(k\pi x_3)$, 
for $i,j,k\in\N_+$, and $x_1,x_2,x_3\in [0,1]$. 
Assuming for example that  $\|Q^{\frac 12}\|_{\LL_2(U,\HH^1_0)}< \infty$, where 
$\HH^1_0:=\HH^1_0(\OO)=\{u\in\HH^1(\OO)\colon u=0\:\:\text{on}\:\:\partial\OO\}$, 
one can get that $\GG Q^{\frac 12}\VV\in D(A)$ for all $\VV\in D(A)$ and 
thus the last inequality in \eqref{con;G} holds.  

For the second example (inspired by \cite{CHJ18}), consider $\GG(\VV)=\VV$ 
for $\VV \in V$, the domain $\OO=[0,1]^3$ and $\epsilon=\mu=1$.
Taking the same orthonormal basis as above, and assuming in addition that $Q^{\frac 12}\in\LL_2(U,\HH^{1+\gamma}(\OO))$ 
with $\gamma>\frac{3}{2}$, one gets for instance
\begin{equation}
\begin{split}
\label{con1;G}
&\|\GG(\VV)\|_{\LL_2(U_0,D(A))}\leq C\|Q^\frac 12\|_{\LL_2(U,\HH^{1+\gamma})}(1+\|\VV\|_{D(A)}).
\end{split}
\end{equation}
Using the definition of the graph norm one gets 
\begin{align*}
\|\GG(\VV)\|_{\LL_2(U_0,D(A))}^2&=
\sum\limits_{k\in\N_+}
\|\VV Q^\frac 12e_k\|_{V}^2+
\sum\limits_{k\in\N_+}
\|A(\VV Q^\frac 12e_k)\|_{V}^2.
\end{align*}
Denoting $\VV=({\bf E}_\VV^T, {\bf H}_\VV^T)^T$ and using the definition 
of the operator $A$, one obtains
\begin{align*}
&\|\GG(\VV)\|_{\LL_2(U_0,D(A))}^2\\
&=
\sum\limits_{k\in \N_+}\sum_{i=1,2,3}\|{\bf E}_\VV^iQ^{\frac 12}e_k\|_U^2
+
\sum\limits_{k\in \N_+}\sum_{i=1,2,3}\|{\bf H}_\VV^iQ^{\frac 12}e_k\|_U^2\\
&\quad+
\sum\limits_{k\in \N_+}\Big(\|\nabla \times ({\bf E}_\VV Q^{\frac 12}e_k)\|_{U^3}^2+\|\nabla \times ({\bf H}_\VV Q^{\frac 12}e_k)\|_{U^3}^2\Big)\\
&\le C\sum\limits_{k\in \N_+}\|Q^{\frac 12}e_k\|_{L^{\infty}(\OO)}^2\|\VV\|_{V}^2
+\sum\limits_{k\in \N_+}\Big(\|\nabla \times ({\bf E}_\VV Q^{\frac 12}e_k)\|_{U^3}^2
+\|\nabla \times ({\bf H}_\VV Q^{\frac 12}e_k)\|_{U^3}^2\Big).
\end{align*}
We now illustrate how to estimate the term $\|\nabla \times ({\bf E}_\VV Q^{\frac 12}e_k)\|_{U^3}^2$ as an example. 
Using the definition of the \textbf{curl} operator, one gets 
\begin{align*}
\|\nabla \times ({\bf E}_\VV Q^{\frac 12}e_k)\|_{U^3}^2
&= 
\|\frac{\partial}{\partial x_2} ({\bf E}_\VV^3Q^{\frac 12}e_k)-\frac{\partial}{\partial x_3} ({\bf E}_\VV^2Q^{\frac 12}e_k) \|_U^2\\
&\quad+
\|\frac{\partial}{\partial x_1} ({\bf E}_\VV^3Q^{\frac 12}e_k)-\frac{\partial}{\partial x_3} ({\bf E}_\VV^1Q^{\frac 12}e_k) \|_U^2\\
&\quad+\|\frac{\partial}{\partial x_1} ({\bf E}_\VV^2Q^{\frac 12}e_k)-\frac{\partial}{\partial x_2} ({\bf E}_\VV^1Q^{\frac 12}e_k) \|_U^2\\
&\le 
C\|Q^{\frac 12}e_k\|^2_{L^{\infty}(\OO)}\Big(\|\frac{\partial}{\partial x_2} {\bf E}_\VV^3-\frac{\partial}{\partial x_3} {\bf E}_\VV^2 \|_U^2
+\|\frac{\partial}{\partial x_1} {\bf E}_\VV^3-\nabla^3 {\bf E}_\VV^1 \|_U^2\\
&\quad
+\|\frac{\partial}{\partial x_1} {\bf E}_\VV^2-\nabla^2 {\bf E}_\VV^1 \|_U^2
\Big)
\\
&\quad+C\Big(\|\frac{\partial}{\partial x_1}Q^{\frac 12}e_k\|_{L^{\infty}(\OO)}^2
+\|\frac{\partial}{\partial x_2}Q^{\frac 12}e_k\|_{L^{\infty}(\OO)}^2+
\|\frac{\partial}{\partial x_3}Q^{\frac 12}e_k\|_{L^{\infty}(\OO)}^2\Big)
\|{\bf E}_\VV \|_{U^3}^2\\
&\le C\|Q^{\frac 12}e_k\|^2_{L^{\infty}(\OO)}
\|\nabla \times {\bf E}_\VV\|_{U^3}^2
+C\|\nabla Q^{\frac 12}e_k\|^2_{L^{\infty}(\OO)}
\| {\bf E}_\VV\|_V^2.
\end{align*}
Combing the above estimates, we obtain 
\begin{align*}
\|\GG(\VV)\|_{\LL_2(U_0,D(A))}^2
&\le 
 C\sum\limits_{k\in \N_+}\|Q^{\frac 12}e_k\|_{L^{\infty}(\OO)}^2
 \bigl(\|\VV\|_{V}^2+\|A \VV\|_{V}^2\bigr)
 +C\sum\limits_{k\in \N_+}\|\nabla Q^{\frac 12}e_k\|_{L^{\infty}(\OO)}^2\|\VV\|_{V}^2.
\end{align*}
Using the Sobolev embedding $\HH^{\gamma}(\OO)\hookrightarrow L^\infty(\OO)$ 
for any $\gamma>\frac{3}{2}$, one finally obtains \eqref{con1;G} and 
the linear growth property of $\GG$.  

The above assumptions suffice to establish well-posedness and 
regularity results of solutions to \eqref{mod;CMAX}. This uses similar arguments as, 
for instance, \cite[Theorem~9]{LSY10} (for a more general drift coefficient in \eqref{mod;CMAX}) 
and \cite[Corollary~3.1]{CHJ18}.
\begin{lm}
\label{lm;ExaBound}
Let $T>0$. Under the Assumptions~\ref{ap;1}-\ref{ap;3}, 
the stochastic Maxwell's equation \eqref{mod;CMAX} is strongly well posed and its 
solution $\UU$ satisfies
\begin{align*}
\EE\left[\sup_{0\leq t\leq T}\|\UU(t)\|_{D(A)}^p\right]<C\left(1+\EE\left[\|\UU(0)\|_{D(A)}^p\right]\right)
\end{align*}
for any~$p\ge 2$. Here, the constant $C$ depends on $p$, $T$, $Q$, bounds for $\FF$ and $\GG$, and $\UU(0)$.
\end{lm}
Subsequently we present a lemma on the H\"older regularity in time of solutions to \eqref{mod;CMAX}. 
This result is important in analysing the approximation error of the proposed time integrator 
in Section~\ref{sect-EXP}. 
\begin{lm}
Let $T>0$. Under the Assumptions~\ref{ap;1}-\ref{ap;3}, the solution $\UU$ of 
the stochastic Maxwell's equation \eqref{mod;CMAX} satisfies 
\label{lm;Holder}
\begin{align*}
\EE\left[\|\UU(t)-\UU(s)\|_V^{2p}\right]\leq C|t-s|^{p},
\end{align*}
for any $0\leq s, t\leq T$, and $p\ge 1$. 
Here, the constant $C$ depends on $p$, $T$, $Q$, bounds for $\FF$ and $\GG$, and $\UU(0)$.
\end{lm}

The proof is very similar to the proof of  \cite[Proposition~3.2]{CHJ18}, 
we omit it for ease of presentation.

Based on the above regularity results for solutions to the stochastic Maxwell's equation \eqref{mod;CMAX}, 
the work \cite{CHJ18} shows mean-square convergence order $\frac12$ of 
the backward Euler--Maruyama scheme (in temporal direction). 
In the next section, we design and analyse 
an explicit and effective numerical scheme, the exponential integrator, 
which has the rate of convergence $1$ and preserves many inherent 
properties of the original problem (in the case of the stochastic Maxwell's equations 
with additive noise).


\section{Exponential integrators for stochastic Maxwell's equations and error analysis}\label{sect-EXP}
This section is concerned with a convergence analysis in 
strong sense of an exponential integrator for the stochastic Maxwell's equation \eqref{mod;CMAX}.
We first show an a priori estimate of the numerical solution. Then the strong convergence rate is studied 
in two cases, first when equation \eqref{mod;CMAX} is driven by additive noise and then for multiplicative 
noise.

Fix a time horizon $T>0$ and an integer $N>0$. Define a stepsize $\Delta t$ such that $T=N\Delta t$. 
We then construct a uniform partition of the interval $[0, T]$ 
\begin{align*}
0=t_0<t_1<\ldots<t_{N-1}<t_N=T
\end{align*} 
with $t_n=n\Delta t$ for $n=0,\ldots,N$. 
Next, we consider the mild solution of the stochastic Maxwell's equation \eqref{mod;CMAX} 
on the small time interval $[t_k,t_{k+1}]$ (with $\UU(t_k)=\UU_k$):
$$
\UU(t_{k+1})=\S(\Delta t)\UU_k+\int_{t_k}^{t_{k+1}}\S(t_{k+1}-s)\FF(\UU(s))\,\text ds
+\int_{t_k}^{t_{k+1}}\S(t_{k+1}-s)\GG(\UU(s))\,\text dW.
$$
By approximating both integrals in the above mild solution at the left end point, one obtains 
the exponential integrator
\begin{align}
\label{sch;exp}
\UU_{k+1}=\S(\Delta t)\UU_k
+\S(\Delta t)\FF(\UU_k)\Delta t
+\S(\Delta t)\GG(\UU_k)\Delta W_k,
\end{align}
where $\Delta W_k=\Delta W(t_{k+1})-\Delta W(t_k)$ stands for Wiener increments.
One readily sees that \eqref{sch;exp} is an explicit numerical approximation of 
the exact solution $\UU(t_{k+1})$ of the stochastic Maxwell's equation \eqref{mod;CMAX}. 

In order to present a result on the strong error of the exponential integrator \eqref{sch;exp}, 
we first show an a priori estimate of the numerical solution.

\begin{tm}
\label{tm1}
Under the Assumptions~\ref{ap;1}-\ref{ap;3}, the numerical solution to the stochastic Maxwell's equation 
given by the exponential integrator \eqref{sch;exp} satisfies 
\begin{align*}
\EE\left[\|\UU_k\|^{2p}_{D(A)}\right]\leq C(\UU_0,Q,T,p,\FF,\GG)
\end{align*}
for all $p\geq1$ and $k=0,1,\ldots,N$. 
\end{tm}
{\bf Proof.} The numerical approximation given by the exponential integrator 
can be rewritten as
\begin{align*}
\UU_k&=\S(t_{k})\UU(0)
+\Delta t\sum\limits_{j=0}^{k-1} 
\S(t_{k}-t_j)\FF(\UU_j)
+\sum\limits_{j=0}^{k-1} 
\S(t_{k}-t_j)\GG(\UU_j)\Delta W_j.
\end{align*}
Taking norm and expectation leads to, for $p\geq1$, 
\begin{align*}
\EE\left[\|\UU_k\|^{2p}_{D(A)}\right]&\leq C\EE\left[\|\S(t_{k})\UU(0)\|^{2p}_{D(A)}\right]
+C\EE\left[\left\|\Delta t\sum\limits_{j=0}^{k-1} 
\S(t_{k}-t_j)\FF(\UU_j)\right\|^{2p}_{D(A)}\right]\\
&\quad+C\EE\left[\left\|\sum\limits_{j=0}^{k-1} 
\S(t_{k}-t_j)\GG(\UU_j)\Delta W_j\right\|^{2p}_{D(A)}\right].
\end{align*}
For the first term, using  the definition of the graph norm and property \eqref{unigro}, we obtain
\begin{align*}
\|\S(t_{k})\UU(0)\|^{2p}_{D(A)}=\left(\|\S(t_{k})\UU(0)\|_V+\|\S(t_{k})A\UU(0)\|_V\right)^{2p}=\|\UU(0)\|^{2p}_{D(A)},
\end{align*}
which leads to $\EE\left[\|\S(t_{k})\UU(0)\|^{2p}_{D(A)}\right]=\EE\left[\|\UU(0)\|^{2p}_{D(A)}\right]$. 
Based on the linear growth property of $\FF$ and H\"{o}lder's inequality, 
the second term is estimated as follows 
\begin{align*}
\left\|\Delta t\sum\limits_{j=0}^{k-1} 
\S(t_{k}-t_j)\FF(\UU_j)\right\|^{2p}_{D(A)}
\leq &C+
C\Delta t^{2p}\left(\sum\limits_{j=0}^{k-1}\|\UU_j\|_{D(A)}\right)^{2p}\\
\leq & C+
C\Delta t^{2p}k^{2p-1}\sum\limits_{j=0}^{k-1}\|\UU_j\|_{D(A)}^{2p}.
\end{align*}
One then obtains
\begin{align*}
\EE\left[\left\|\Delta t\sum\limits_{j=0}^{k-1} 
\S(t_{k}-t_j)\FF(\UU_j)\right\|^{2p}_{D(A)}\right]
\leq C+
C\Delta t\EE\left[\sum\limits_{j=0}^{k-1}\|\UU_j\|_{D(A)}^{2p}\right].
\end{align*}
The third term is equivalent to 
\begin{align*}
\EE\left[\left\|\sum\limits_{j=0}^{k-1} 
\S(t_{k}-t_j)\GG(\UU_j)\Delta W_j\right\|^{2p}_{D(A)}\right]
&=\EE\left[\left\|\int_{0}^{t_{k}}
\S\left(t_{k}-[\frac s{\Delta t}]\Delta t\right)\GG(\UU_{[\frac s{\Delta t}]\Delta t})\,\text d W(s)\right\|^{2p}_{D(A)}\right]
\end{align*}
with $[\frac s{\Delta t}]$ being the integer part of $\frac s{\Delta t}$. 
The Burkholder--Davis--Gundy inequality for stochastic integrals and our assumption on $\GG$ give
\begin{align*}
&\EE\left[\left\|\int_{0}^{t_{k}}
\S\left(t_{k}-[\frac s{\Delta t}]\Delta t\right)\GG(\UU_{[\frac s{\Delta t}]\Delta t})\,\text dW(s)\right\|^{2p}_{D(A)}\right]\leq \\
&\leq C\EE\left[\left(\int_{0}^{t_{k}}
\left\|\GG(\UU_{[\frac s{\Delta t}]\Delta t})\right\|_{\LL_2(U_0,D(A))}^2\,\text ds\right)^{p}\right]\\
&\leq C+C\EE\left[\left(\int_{0}^{t_{k}}
\left\|\UU_{[\frac s{\Delta t}]\Delta t}\right\|_{D(A)}^2\,\text ds\right)^{p}\right]
=C+C\EE\left[\left(\Delta t\sum\limits_{j=0}^{k-1} \|\UU_j\|_{D(A)}^2\right)^{p}\right].
\end{align*}
Using H\"{o}lder's inequality, the last term in the above inequality becomes
\begin{align*}
\left(\Delta t\sum\limits_{j=0}^{k-1} \|\UU_j\|_{D(A)}^2\right)^{p}\leq 
\Delta t^{p}k^{p-1}\sum\limits_{j=0}^{k-1} \|\UU_j\|_{D(A)}^{2p}.
\end{align*}
Taking expectation, we then obtain
\begin{align*}
\EE\left[\left\|\int_{0}^{t_{k}}
\S\left(t_{k}-[\frac s{\Delta t}]\Delta t\right)\GG(\UU(s))\,\text dW(s)\right\|^{2p}_{D(A)}\right]
\leq C+ C\Delta t\sum\limits_{j=0}^{k-1} \EE\left[\|\UU_j\|_{D(A)}^{2p}\right].
\end{align*}
Altogether, we get that
\begin{align*}
\EE\left[\|\UU_k\|^{2p}_{D(A)}\right]\leq C+C\Delta t\EE\left[\sum\limits_{j=0}^{k-1}\|\UU_j\|_{D(A)}^{2p}\right].
\end{align*} 
A discrete Gronwall inequality concludes the proof.\hfill$\square$\\
Using the above theorem, we arrive at
\begin{cor}\label{cor:sup}
Under the same assumptions as in Theorem~\ref{tm1}, for all $p\geq1$, 
there exists a constant $C:=C(\UU(0),Q,T,p,\FF,\GG)$ such that
\begin{align}
\label{Nmm}
\EE\left[\sup_{0\leq k\leq N}\|\UU_k\|^{2p}_{D(A)}\right]\leq C.
\end{align}
\end{cor}
{\bf Proof.}
The main idea to derive the estimate \eqref{Nmm} is to properly estimate the stochastic integral
\begin{align*}
&\EE\left[\sup_{0\leq k\leq N}\left\|\sum\limits_{j=0}^{k-1} 
\S(t_{k}-t_j)\GG(\UU_j)\Delta W_j\right\|^{2p}_{D(A)}\right]=\\
&=\EE\left[\sup_{0\leq k\leq N}
\left\|
\int_{0}^{t_{k}}
\S\left(t_{k}-[\frac s{\Delta t}]\Delta t\right)\GG(\UU_{[\frac s{\Delta t}]\Delta t})
\,\text dW(s)\right\|^{2p}_{D(A)}\right].
\end{align*}
Based on the unitarity of $S(\cdot)$, Burkholder--Davis--Gundy's inequality, 
H\"older's inequality, and our assumptions on $\GG$,  
the right hand side (RHS) of the above equality becomes
\begin{align*}
{\rm RHS}&\leq C\EE\left[\left(\int_{0}^{T}
\left\|\GG(\UU_{[\frac s{\Delta t}]\Delta t})\right\|_{\LL_2(U_0,D(A))}^2\,\text ds\right)^{p}\right]\\
& \leq
C+C\Delta t\sum\limits_{j=0}^{N-1} \EE\left[\|\UU_j\|_{D(A)}^{2p}\right]\leq C,
\end{align*}
where we use the result of Theorem~\ref{tm1} in the last step. 
The estimations of the other terms in the numerical solution are done in a similar way as in the previous result. 
\hfill$\square$

We are now in position to show the error estimates of the exponential integrator for the stochastic Maxwell's equation 
\eqref{mod;CMAX} driven by additive noise.

\begin{tm}\label{thm-strong-add}
Let Assumptions~\ref{ap;1}-\ref{ap;3} hold. Assume in addition that $\FF\in C_b^2(V)$ and 
$\GG$ does not dependent on $\UU$. 
The strong error of the exponential integrator \eqref{sch;exp} when applied 
to the stochastic Maxwell's equation \eqref{mod;CMAX} verifies, for all $p\geq1$, 
\begin{align*}
\left(\EE\left[ \max_{k=0,\ldots,N}\|\UU(t_{k})-\UU_{k}\|_V^{2p}\right]
\right)^\frac{1}{2p}
\leq C\Delta t,
\end{align*} 
where the positive constant $C$ depends on bounds for $\FF$ (and its derivatives) and $\GG$, 
as well as on $T$, $p$ and $Q$.
\end{tm}


{\bf Proof.}
Let us denote $\epsilon_{k}=\UU(t_{k})-\UU_{k}$, for $k=0,\ldots,N$. We then have
\begin{align}
\label{err}
\epsilon_{k+1}&=
\sum\limits_{j=0}^k 
\int_{t_j}^{t_{j+1}}
\left(\S(t_{k+1}-s)\FF(\UU(s))
-\S(t_{k+1}-t_j)\FF(\UU_j)\right)\,\text ds\nonumber\\
&\quad+\sum\limits_{j=0}^k 
\int_{t_j}^{t_{j+1}}
\left((\S(t_{k+1}-s)
-\S(t_{k+1}-t_j))\GG\right) \,\text dW(s)\nonumber\\
&=:Err_1^k+Err_2^k.
\end{align}
We now rewrite the term $Err_1^k$ as
\begin{align*}
Err_1^k&=\sum\limits_{j=0}^k 
\int_{t_j}^{t_{j+1}}
\left(\S(t_{k+1}-s)(\FF(\UU(s))-\FF(\UU(t_j)))\right)\,\text ds\\
&\quad+\sum\limits_{j=0}^k 
\int_{t_j}^{t_{j+1}}
\left(\left(\S(t_{k+1}-s)-\S(t_{k+1}-t_j)\right)
\FF(\UU(t_j))\right)\,\text ds\\
&\quad+\sum\limits_{j=0}^k 
\int_{t_j}^{t_{j+1}}
\left(\S(t_{k+1}-t_j)(\FF(\UU(t_j))-\FF(\UU_j))\right)\,\text ds\\
&=:\MyRoman{1}_1^k+\MyRoman{1}_2^k+\MyRoman{1}_3^k.
\end{align*}
We first estimate the term $\MyRoman{1}_1^k$. 
Using a Taylor expansion, we obtain 
\begin{align*}
\FF(\UU(s))-\FF(\UU(t_j))&=\frac{\partial \FF}{\partial u}(\UU(t_j))(\UU(s)-\UU(t_j))\\
&\quad+\frac{1}{2}\frac{\partial^2 \FF}{\partial u^2}(\Theta)
(\UU(s)-\UU(t_j),\UU(s)-\UU(t_j)),
\end{align*}
where $\Theta:=\theta\UU(s)+(1-\theta)\UU(t_j)$, for some $\theta\in[0,1]$, depends on $\UU(s)$ and $\UU(t_j)$. 
Combing this with the mild formulation of the exact solution on the interval $[t_j,s]$,  
\begin{align*}
\UU(s)=\S(s-t_j)\UU(t_j)+\int_{t_j}^s\S(s-r)\FF(\UU(r))\,\text dr+
\int_{t_j}^s\S(s-r)\GG\,\text dW(r),
\end{align*}
we rewrite the term $\MyRoman{1}_1^k$ as
$$
{\MyRoman{1}_1^k}=\AA_1^k+\AA_2^k,
$$
where we define 
\begin{align*}
\AA_1^k&=\sum\limits_{j=0}^k 
\int_{t_j}^{t_{j+1}}
\S(t_{k+1}-s)\frac{\partial \FF}{\partial u}(\UU(t_j))(\S(s-t_j)-Id)\UU(t_j)\,\text ds\\
&\quad+\sum\limits_{j=0}^k 
\int_{t_j}^{t_{j+1}}
\S(t_{k+1}-s)\frac{\partial \FF}{\partial u}(\UU(t_j))\int_{t_j}^s\S(s-r)\FF(\UU(r))\,\text dr\,\text ds\\
&\quad+\sum\limits_{j=0}^k 
\int_{t_j}^{t_{j+1}}
\S(t_{k+1}-s)\frac{\partial \FF}{\partial u}(\UU(t_j))\int_{t_j}^s\S(s-r)\GG\,\text dW(r)\,\text ds\\
&=:\MyRoman{2}_1^k+\MyRoman{2}_2^k+\MyRoman{2}_3^k,
\end{align*}
and
\begin{align*}
\AA_2^k=\sum\limits_{j=0}^k 
\int_{t_j}^{t_{j+1}}
\S(t_{k+1}-s)\frac{1}{2}\frac{\partial^2 \FF}{\partial u^2}(\Theta)
(\UU(s)-\UU(t_j),\UU(s)-\UU(t_j))\,\text ds.
\end{align*}
The assumption that $\FF\in C_b^2(V)$ and the H\"older continuity of the exact solution $\UU$ 
in Lemma~\ref{lm;Holder} provide us with the bound $\EE\left[\|\AA_2\|_V^{2p}\right]\le C\Delta t^{2p}$. 
For the term $\MyRoman{2}_1$, we use property \eqref{unigro}, the boundedness of the derivatives 
of $\FF$ and Lemma~\ref{lm;SG}, combined with H\"{o}lder's inequality, to deduce that 
\begin{align*}
{\|\MyRoman{2}_1^k\|}_V&\leq \sum\limits_{j=0}^k 
\int_{t_j}^{t_{j+1}}
\left\|\frac{\partial \FF}{\partial u}(\UU(t_j))(\S(s-t_j)-Id)\UU(t_j)\right\|_V\,\text ds\\
&\leq C\sum\limits_{j=0}^k 
\int_{t_j}^{t_{j+1}}|s-t_j|
\|\UU(t_j)\|_{D(A)}\,\text ds
\leq  C(\Delta t)^2 \sum\limits_{j=0}^k \|\UU(t_j)\|_{D(A)}\\
&\leq C(\Delta t)^2 \left(\sum\limits_{j=0}^k \|\UU(t_j)\|_{D(A)}^{2p}\right)^\frac{1}{2p}
\left(\frac{t_{k+1}}{\Delta t}\right)^{\frac{2p-1}{2p}}\\
&\leq  C\Delta t
\left(\sup_{0\leq j\leq k}\|\UU(t_j)\|_{D(A)}^{2p}\right)^\frac{1}{2p}.
\end{align*}
This leads to
\begin{align*}
\EE\left[\max_{k=0,\ldots,N-1}{\|\MyRoman{2}_1^k\|}_V^{2p}\right]\leq
C(\Delta t)^{2p}\EE\left[\sup_{0\leq j\leq N }\left\|\UU(t_j)\right\|_{D(A)}^{2p}\right]
\leq C(\Delta t)^{2p}
\end{align*}
using Lemma~\ref{lm;ExaBound}.
Next, we estimate the term $\MyRoman{2}_2^k$. 
Using Lemma~\ref{lm;SG} and H\"{o}lder's inequality, we obtain
\begin{align*}
{\|\MyRoman{2}_2^k\|}_V&\leq
C\sum\limits_{j=0}^k 
\int_{t_j}^{t_{j+1}}
\int_{t_j}^s\|\FF(\UU(r))\|_V\,\text dr\,\text ds\\
&\leq C\sum\limits_{j=0}^k 
\int_{t_j}^{t_{j+1}}
\int_{t_j}^s(1+\|\UU(r)\|_V)\,\text dr\,\text ds\\
&\leq C\Delta t+
C\sum\limits_{j=0}^k 
\int_{t_j}^{t_{j+1}} (s-t_j)^{\frac{2p-1}{2p}}
\left(\int_{t_j}^s\|\UU(r)\|_V^{2p}\,\text dr\right)^\frac {1}{2p}\,\text ds\\
&\leq 
C\Delta t+
C\Delta t
\left(\sup_{0\leq t\leq T}\|\UU(t)\|_V^{2p}\right)^\frac {1}{2p}.
\end{align*}
From Lemma~\ref{lm;ExaBound}, It then follows that 
\begin{align*}
\EE\left[\max_{k=0,\ldots,N-1}{\|\MyRoman{2}_2^k\|}_V^{2p}\right]\leq
C(\Delta t)^{2p}
+C(\Delta t)^{2p}\EE\left[\sup_{0\leq t\leq T}\left\|\UU(t)\right\|_V^{2p}\right]
\leq C(\Delta t)^{2p}.
\end{align*}
We now proceed to the estimation of the term $\MyRoman{2}_3^k$. 
First notice that stochastic Fubini's theorem leads to
\begin{align*}
\MyRoman{2}_3^k=
&\sum\limits_{j=0}^k 
\int_{t_j}^{t_{j+1}}
\S(t_{k+1}-s)\frac{\partial \FF}{\partial u}(\UU(t_j))\int_{t_j}^s\S(s-r)\GG\,\text dW(r)ds\\
&=\sum\limits_{j=0}^k 
\int_{t_j}^{t_{j+1}}\int_r^{t_{j+1}}
\S(t_{k+1}-s)\frac{\partial \FF}{\partial u}(\UU(t_j))
\S(s-r)\,\text ds\,\text dW(r)\\
&=\int_0^{t_{k+1}}\int_r^{([\frac r{\Delta t}]+1)\Delta t}\S(t_{k+1}-s)\frac{\partial \FF}{\partial u}(\UU([\frac s{\Delta t}]\Delta t))
\S(s-r)\,\text ds\,\text dW(r)
\end{align*}
and the integrand in the above equation is $\mathcal F_r$-adaptive. 
Then by the Burkholder--Davis--Gundy's inequality, we get
\begin{align*}
&\E[\max_{k=0,\ldots,N-1}\|\MyRoman{2}_3^k\|_V^{2p}]\\
&\le 
C\EE\left[
\left(\int_0^T\left\|\int_r^{([\frac r{\Delta t}]+1)\Delta t}\S(t_{k+1}-s)\frac{\partial \FF}{\partial u}(\UU({[\frac s{\Delta t}]\Delta t}))
\S(s-r)\,\text ds\right\|_{\LL_2(U_0,V)}^2\,\text dr\right)^{p}\right].
\end{align*}
Then, using the assumption that $\FF\in C^2_b(V)$, we obtain 
\begin{align*}
&\EE\left[\max_{k=0,\ldots,N-1}{\|\MyRoman{2}_3^k\|}_V^{2p}\right]
\nonumber\\
\leq &
C\EE\left[
\left(\sum\limits_{j=0}^{N-1}\int_{t_j}^{t_{j+1}}
\Bigg{(}\int_r^{t_{j+1}}
\left\|\S(t_{k+1}-s)\frac{\partial \FF}{\partial u}
(\UU(t_j))
\S(s-r)\right\|_{\LL_2(U_0,V)}\,\text ds\Bigg{)}^2\,\text dr\right)^{p}\right]
\nonumber\\
\leq &
C\EE\left[
\left(\sum\limits_{j=0}^{N-1}\int_{t_j}^{t_{j+1}}
\Bigg{(}\int_r^{t_{j+1}}
\left\|Q^\frac 12\right\|_{\LL_2(U,V)}\,\text ds\Bigg{)}^2\,\text 
dr\right)^{p}\right]
\leq 
C(\Delta t)^{2p}.
\end{align*}
Thus, the above allows us to get the following estimate
\begin{align*}
\EE\left[\max_{k=0,\ldots,N-1}{\|\AA_1\|}_V^{2p}\right]\leq C (\Delta t)^{2p},
\end{align*}
which implies the estimate
\begin{align*}
\EE\left[\max_{k=0,\ldots,N-1}{\|\MyRoman{1}_1^k\|}_V^{2p}\right]\leq C (\Delta t)^{2p}.
\end{align*}
For the term $\MyRoman{1}_2^k$, we use the unitary property of 
the semigroup \eqref{unigro} to get
\begin{align*}
{\|\MyRoman{1}_2^k\|}_V&\leq
\sum\limits_{j=0}^k 
\int_{t_j}^{t_{j+1}}
\left\|(\S(t_{k+1}-s)-\S(t_{k+1}-t_j))
\FF(\UU(t_j))\right\|_V\,\text ds\\
&=\sum\limits_{j=0}^k 
\int_{t_j}^{t_{j+1}}
\left\|(\S(t_{j}-s)-Id)
\FF(\UU(t_j))\right\|_V\,\text ds.
\end{align*}
According to Lemma~\ref{lm;SG} and the linear growth property of $\FF$, 
the above term can be bounded by
\begin{align*}
{\|\MyRoman{1}_2^k\|}_V
&\leq C \sum\limits_{j=0}^k 
\int_{t_j}^{t_{j+1}}
|t_{j}-s|
\left\|\FF(\UU(t_j))\right\|_{D(A)}\,\text ds\\
&\leq C (\Delta t)^2 \sum\limits_{j=0}^k
\left\|\FF(\UU(t_j))\right\|_{D(A)}\\
&\leq C\Delta t
+C (\Delta t)^2 \sum\limits_{j=0}^k
\left\|\UU(t_j)\right\|_{D(A)}.
\end{align*}
Taking the $2p$-th power on both sides of the above inequality and then expectation, we obtain 
\begin{align*}
\EE\left[\max_{k=0,\ldots,N-1}{\|\MyRoman{1}_2^k\|}_V^{2p}\right]\leq
C(\Delta t)^{2p}
+C(\Delta t)^{2p}\EE\left[\sup_{0\leq t\leq T}\left\|\UU(t)\right\|_{D(A)}^{2p}\right]
\leq C(\Delta t)^{2p}
\end{align*}
by Lemma~\ref{lm;ExaBound} in Section~\ref{sec;2}.
For the term $\MyRoman{1}_3^k$, similarly as above, using properties of the semigroup and of $\FF$, and 
H\"{o}lder's inequality, we obtain
\begin{align*}
{\|\MyRoman{1}_3^k\|}_V
&\leq \Delta t\sum\limits_{j=0}^k 
\|\epsilon_j\|_V
\leq \Delta t
\left(\sum\limits_{j=0}^k 
\|\epsilon_j\|_V^{2p}\right)^\frac 1{2p}
\left(\frac{t_{k+1}}{\Delta t}\right)^{\frac{2p-1}{2p}}\\
&\leq C \Delta t
\left(\sum\limits_{j=0}^k 
\|\epsilon_j\|_V^{2p}\right)^\frac 1{2p}
(\Delta t)^{\frac{1-2p}{2p}}
=C(\Delta t)^{\frac{1}{2p}} 
\left(\sum\limits_{j=0}^k 
\|\epsilon_j\|_V^{2p}\right)^\frac 1{2p}.
\end{align*}
This gives us
\begin{align*}
\EE\left[\max_{k=0,\ldots,N-1}{\|\MyRoman{1}_3^k\|}_V^{2p}\right]\leq
C\Delta t\sum\limits_{j=0}^{N-1} 
\EE\left[\max_{l=0,\ldots,j}\|\epsilon_l\|_V^{2p}\right].
\end{align*}
The last term $Err_2^k$ can be bounded as follows
\begin{align*}
&\EE\left[\max_{k=0,\ldots,N-1}\|Err_2^k\|_V^{2p}\right]\\
&=
\EE\left[\max_{k=0,\ldots,N-1}
\left\|
\sum\limits_{j=0}^k 
\int_{t_j}^{t_{j+1}}
(\S(t_{k+1}-s)-\S(t_{k+1}-t_j))\GG\,\text d W(s)\right\|_V^{2p}\right]\\
&=\EE\left[\max_{k=0,\ldots,N-1}
\left\|
\int_{0}^{t_{k+1}}
\left(\S(t_{k+1}-s)-\S(t_{k+1}-\left[\frac{s}{\Delta t}\right]\Delta t)\right)\GG\,\text dW(s)\right\|_V^{2p}\right]\\
&\leq \EE\left[\sup_{0\leq t\leq T}
\left\|
\int_{0}^{t}
\S(t-\left[\frac{s}{\Delta t}\right]\Delta t)
\left(\S(\left[\frac{s}{\Delta t}\right]\Delta t-s)-Id\right)\GG\,\text dW(s)\right\|_V^{2p}\right].
\end{align*}
Thanks to Burkholder--Davis--Gundy's inequality and properties of the semigroup, we obtain
\begin{align*}
\EE\left[\max_{k=0,\ldots,N-1}\|Err_2^k\|_V^{2p}\right]
&\leq  C\EE\left[\left(
\int_{0}^{T}
\left\|
(\S(\left[\frac{s}{\Delta t}\right]\Delta t-s)-Id)\GG \right\|_{\LL_2(U_0,V)}^2
\,\text ds\right)^{p}\right]\\
&= 
C\EE\left[\left(\sum\limits_{j=0}^{N-1}\int_{t_j}^{t_{j+1}}
\left\|(\S(t_j-s)-Id)\GG \right\|_{\LL_2(U_0,V)}^2
\,\text ds\right)^{p}\right]\\
&\leq  C\EE\left[\left(\sum\limits_{j=0}^{N-1}\int_{t_j}^{t_{j+1}}
|t_j-s|^2\left\|\GG Q^\frac 12\right\|_{\LL_2(U,D(A))}^2
\,\text ds\right)^{p}\right]\\
&\leq C(\Delta t)^{2p},
\end{align*}
where we have used the linear growth property of $\GG$ in $\LL_2(U_0,D(A))$ in the last step.

Collecting all the above estimates gives us the bound
\begin{align*}
\EE\left[\max_{k=0,\ldots,N-1}\|\epsilon_{k+1}\|_V^{2p}\right]\leq C(\Delta t)^{2p}+
C\Delta t\sum\limits_{j=0}^{N-1} 
\EE\left[\max_{l=0,\ldots,j}\|\epsilon_l\|_V^{2p}\right].
\end{align*}
An application of Gronwall's inequality yields
\begin{align*}
\left(\EE\left[\max_{k=0,\ldots,N}\|\epsilon_{k}\|_V^{2p}\right]
\right)^\frac{1}{2p}
\leq C\Delta t,
\end{align*}
which means that the strong order of the exponential scheme is $1$ if the noise is additive in the stochastic Maxwell's equation \eqref{mod;CMAX}.
\hfill$\square$

Now we turn to the case where the stochastic Maxwell's equation \eqref{mod;CMAX} 
is driven by a more general multiplicative noise.

\begin{tm}\label{thm-strong}
Let Assumptions \ref{ap;1}-\ref{ap;3} hold. 
The strong error of the exponential integrator \eqref{sch;exp} when applied 
to the stochastic Maxwell's equation \eqref{mod;CMAX} verifies, for all $p\geq1$, 
\begin{align*}
\left(\EE \left[\max_{k=0,\ldots,N}\|\UU(t_{k})-\UU_{k}\|_V^{2p}\right]
\right)^\frac{1}{2p}
\leq C\Delta t^\frac 12,
\end{align*} 
where the positive constant $C$ depends on the Lipschitz coefficients of $\FF$ and $\GG$, 
$p$, $\UU(0)$, $Q$ and $T$. 
\end{tm}
{\bf Proof.}
When the noise is multiplicative, the term $Err_2^k$ in \eqref{err} becomes 
\begin{align*}
Err_2^k=\sum\limits_{j=0}^k 
\int_{t_j}^{t_{j+1}}
\left(\S(t_{k+1}-s)\GG(\UU(s))
-\S(t_{k+1}-t_j)\GG(\UU_j)\right)\,\text dW(s),
\end{align*}
which can be rewritten as
\begin{align*}
Err_2^k&=\sum\limits_{j=0}^k 
\int_{t_j}^{t_{j+1}}
\S(t_{k+1}-s)(\GG(\UU(s))-\GG(\UU(t_j)))\,\text dW(s)\\
&\quad+\sum\limits_{j=0}^k 
\int_{t_j}^{t_{j+1}}
\left(\S(t_{k+1}-s)-\S(t_{k+1}-t_j)\right)
\GG(\UU(t_j))
\,\text dW(s)\\
&\quad+\sum\limits_{j=0}^k 
\int_{t_j}^{t_{j+1}}
\S(t_{k+1}-t_j)(\GG(\UU(t_j))-\GG(\UU_j))\,\text dW(s)\\
&=:\MyRoman{3}_1+\MyRoman{3}_2+\MyRoman{3}_3.
\end{align*}
By Burkholder--Davis--Gundy's inequality and the assumptions on $\GG$, one obtains
\begin{align*}
&\EE \left[\max_{k=0,\ldots,N-1}\|\MyRoman{3}_1\|_V^{2p}\right]\\\leq&
\EE \left[\sup_{0\leq t\leq T}
\left\| 
\int_{0}^{t}
\S(t-s)(\GG(\UU(s))-\GG(\UU(\left[\frac{s}{\Delta t}\right]\Delta t)))\,\text dW(s)
\right\|_V^{2p}\right]\\
\leq&
C\EE\left[\left(
\int_{0}^{T}
\left\|
\GG(\UU(s))-\GG(\UU(\left[\frac{s}{\Delta t}\right]\Delta t))
\right\|_{\LL_2(U_0,V)}^2
\,\text ds\right)^{p}\right]\\
\leq & C\EE\left[\left(
\int_{0}^{T}
\left\|
\UU(s)-\UU(\left[\frac{s}{\Delta t}\right]\Delta t)
\right\|_V^2
\,\text ds\right)^{p}\right].
\end{align*}
Based on H\"{o}lder's inequality and the continuity of 
$\UU$ in Lemma~\ref{lm;Holder}, we have
\begin{align*}
&\EE \left[\max_{k=0,\ldots,N-1}\|\MyRoman{3}_1\|_V^{2p}\right]\\
\leq& C\EE\left[\left(
\left(\int_{0}^{T}
\left\|
\UU(s)-\UU(\left[\frac{s}{\Delta t}\right]\Delta t)
\right\|_V^{2p}
\,\text ds
\right)^\frac 1p
T^\frac{p-1}{p}
\right)^{p}\right]\\
\leq&C \EE\left[\int_{0}^{T}
\left\|
\UU(s)-\UU(\left[\frac{s}{\Delta t}\right]\Delta t)
\right\|_V^{2p}
\,\text ds\right]\\
\leq&
C\sum\limits_{j=0}^{N-1} 
\int_{t_j}^{t_{j+1}}|s-t_j|^{p}
\,\text ds\leq C(\Delta t)^p.
\end{align*}
Similarly, for the term $\MyRoman{3}_2$, we obtain
\begin{align*}
&\EE \left[{\max_{k=0,\cdots,N-1}\|\MyRoman{3}_2\|_V^{2p}}\right]\\
\leq &
\EE\left[ \sup_{0\leq t\leq T}
\left(\left\| 
\int_{0}^{t}
\S(t-s)-\S(t-\left[\frac{s}{\Delta t}\right]\Delta t)\right)
\GG(\UU(\left[\frac{s}{\Delta t}\right]\Delta t))
\,\text dW(s)
\right\|_V^{2p}\right]\\
\leq &
C\EE\left[ \left(\int_{0}^{T}
\left\|
\left(\S(t-s)-\S(t-\left[\frac{s}{\Delta t}\right]\Delta t)\right)
\GG(\UU(\left[\frac{s}{\Delta t}\right]\Delta t))
\right\|_{\LL_2(U_0,V)}^2
\,\text ds\right)^p\right]\\
\leq&
CT^{p-1}\EE \left[\int_{0}^{T}
\left\|(\S(s-\left[\frac{s}{\Delta t}\right]\Delta t)-Id)
\GG(\UU(\left[\frac{s}{\Delta t}\right]\Delta t))
\right\|_{\LL_2(U_0,V)}^{2p}
\,\text ds\right]\\
\leq& C\sum\limits_{j=0}^{N-1} 
\int_{t_j}^{t_{j+1}} 
|s-\left[\frac{s}{\Delta t}\right]\Delta t|^{2p}
\EE\left[\left\|\GG(\UU(\left[\frac{s}{\Delta t}\right]\Delta t))
\right\|_{\LL_2(U_0,D(A))}^{2p}
\,\text ds\right]\\
\leq& C (\Delta t)^{2p}.
\end{align*}
For the last term $\MyRoman{3}_3$, using Assumption \ref{ap;3}, we get
\begin{align*}
&\EE\left[ {\max_{k=0,\ldots,N-1}\|\MyRoman{3}_3\|_V^{2p}}\right]\\
\leq&
\EE\left[ \sup_{0\leq t\leq T}
\left\| 
\int_{0}^{t} 
\S(\left[\frac{s}{\Delta t}\right]\Delta t)(\GG(\UU(\left[\frac{s}{\Delta t}\right]\Delta t))-\GG(\UU_{\left[\frac{s}{\Delta t}\right]}))\,\text dW(s)
\right\|_V^{2p}\right]\\
\leq &C
\EE \left[\left(
\int_{0}^{T} 
\left\|
\GG(\UU(\left[\frac{s}{\Delta t}\right]\Delta t))-\GG(\UU_{\left[\frac{s}{\Delta t}\right]})
\right\|^2_{\LL_2(U_0,V)}\,\text ds
\right)^p\right]\\
\leq &
C \EE\left[ \left(
\int_{0}^{T} 
\left\|
\UU(\left[\frac{s}{\Delta t}\right]\Delta t)-\UU_{\left[\frac{s}{\Delta t}\right]}
\right\|^2_V\,\text ds
\right)^p\right]\\
\leq &C\Delta t 
\sum\limits_{j=0}^{N-1} \EE\left[\max_{l=0,\ldots,j}
\left\|
\UU(t_l)-\UU_l
\right\|^{2p}_V\right].
\end{align*}
Altogether, we obtain
\begin{align*}
\EE\left[\max_{k=0,\ldots,N-1}{\|Err_2^k\|}_V^{2p}\right]\leq
C(\Delta t)^{p}
+C\Delta t\sum\limits_{j=0}^{N-1} 
\EE\left[\max_{l=0,\ldots,j}\|\epsilon_l\|_V^{2p}\right],
\end{align*}
where we recall the notation $\epsilon_l=\UU(t_l)-\UU_l$.
Another difference with the proof for the additive noise case is estimating the term $\MyRoman{1}_1^k$.
Using \eqref{unigro} and Assumption~\ref{ap;2}, we obtain
\begin{align*}
{\|\MyRoman{1}_1^k\|}_V^{2p}&\leq
\left(\sum\limits_{j=0}^k 
\int_{t_j}^{t_{j+1}}
\| \S(t_{k+1}-s)(\FF(\UU(s))-\FF(\UU(t_j)))\|_V\,\text ds\right)^{2p}\\
&\leq C
\left(\sum\limits_{j=0}^k 
\int_{t_j}^{t_{j+1}}
\|\UU(s)-\UU(t_j)\|_V\,\text ds\right)^{2p}\\
&\leq C
\sum\limits_{j=0}^k 
\int_{t_j}^{t_{j+1}}
\|\UU(s)-\UU(t_j)\|_V^{2p}\,\text ds.
\end{align*}
Using Lemma~\ref{lm;Holder}, one gets
\begin{align*}
\EE\left[\max_{k=0,\ldots,N-1}{\|\MyRoman{1}_1^k\|}_V^{2p}\right]
\leq&C
\sum\limits_{j=0}^k 
\int_{t_j}^{t_{j+1}}
|s-t_j|^{p}\,\text ds\leq
 C (\Delta t)^{p}.
 \end{align*}
Putting all these estimates together yields
\begin{align*}
\EE \left[\max_{k=0,\ldots,N-1}\|\epsilon_{k+1}\|_V^{2p}\right]\leq C(\Delta t)^{p}+
C\Delta t\sum\limits_{j=0}^{N-1} 
\EE\left[\max_{l=0,\ldots,j}\|\epsilon_l\|_V^{2p}\right].
\end{align*}
An application of Gronwall's inequality completes the proof, that is, on gets
\begin{align*}
\left(\EE \left[\max_{k=0,\ldots,N}\|\epsilon_{k}\|_V^{2p}\right]
\right)^\frac{1}{2p}
\leq C(\Delta t)^\frac 12.
\end{align*}
\hfill$\square$

\section{Linear stochastic Maxwell's equations with additive noise}\label{sect-LSM}
In this section, we study phenomena where the densities of the electric 
and magnetic currents are assumed to be linear. 
This is an important example of application of stochastic Maxwell's equations in physics, see e.g. \cite[Chapter 3, pages 112-114]{RKT89}. 
We thus now inspect the long-time behavior of the exponential integrator applied 
to the linear stochastic Maxwell's equation with additive noise. 
We also briefly comment on the symplectic structure of the exact and numerical solutions. 
For simplicity of presentation, in this section we consider a similar setting as in \cite{CHZ16}:  
we assume that $\epsilon=\mu=1$, take $\FF=0$ and 
$\GG=(\lambda_1,\lambda_1,\lambda_1,\lambda_2,\lambda_2,\lambda_2)^\top$ for two real numbers $\lambda_1$ and $\lambda_2$.  
Then the stochastic Maxwell's equation \eqref{mod;CMAX} becomes the linear stochastic Maxwell's equation with additive noise: 
\begin{align}
\label{mod;Max1}
\text d\E-\nabla\times\H \,\text dt&=\lambda_1{\bf e}\,\text dW,\nonumber\\
\text d\H+\nabla\times\E \,\text dt&=\lambda_2{\bf e}\,\text dW,
\end{align}
where  ${\bf{e}}=(1,1,1)^\top$.
In \cite{CHZ16}, it is shown that the averaged energy increases linearly with respect to the evolution of time and 
that the flow of the linear stochastic Maxwell's equation with additive noise preserves the divergence in the sense of expectation. 
We now recall these results and analyse the behavior of the exponential integrator with respect to the preservation of these 
geometric properties of the problem. 
\begin{lm}[Theorems~2.1~and~2.2 in \cite{CHZ16}, Theorem~3.1 in \cite{CHJ18a}]
Consider the linear stochastic Maxwell's equation \eqref{mod;Max1} with a trace class noise. 
There exists a constant $K=3\left(\lambda_1^2+\lambda_2^2\right)\Tr(Q)$ such that the averaged 
energy of the exact solution satisfies the trace formula
\begin{align*}
\EE\left[\Phi^{exact}(t)\right]=\EE\left[\Phi^{exact}(0)\right]+Kt\quad\text{for all times}\quad t,
\end{align*}
where $\displaystyle\Phi^{exact}(t):=\int_{\OO}\left(\|\E(t)\|^2+\|\H(t)\|^2\right)\,\text d{\bf x}$ 
denotes the energy of the problem. 

Assume that $Q^\frac 12\in \LL(\LL^2(\OO),\HH^1(\OO))$, then the solution to 
equation \eqref{mod;Max1} preserves the averaged divergence 
\begin{align*}
\EE\left[\ddiv (\E(t))\right]=\EE\left[\ddiv (\E(0))\right],\quad 
\EE\left[\ddiv (\H(t))\right]=\EE\left[\ddiv (\H(0))\right]
\end{align*}
for all times $t$.

The solutions to Maxwell's equation \eqref{mod;Max1} preserves the symplectic structure
$$
\overline\omega(t)=\overline\omega(0)\quad\mathbb{P}\text{-a.s.},
$$
where $\overline\omega(t):=\displaystyle\int_{\OO}\text d\E(t,{\bf x})\wedge \text d\H(t,{\bf x})\,\text d\bf x$.
\end{lm}

We now show that the proposed exponential integrator possesses the same long-time behavior 
as the exact solution to the linear stochastic Maxwell's equation. 
This is certainly not the case for traditional time integrators such as 
Euler--Maruyama's scheme, see the numerical experiments below.  
Recall, that under this setting, the exponential integrator applied to \eqref{mod;Max1} reads
\begin{align}
\label{sm1}
\UU_{k+1}=\S(\Delta t)\UU_k
+\S(\Delta t)\GG\Delta W_k.
\end{align}
We look at the trace formula for the energy first. 
\begin{prop}
\label{tm;AvEn}
The numerical scheme \eqref{sm1} satisfies the same trace formula for the energy as the exact solution 
to the linear stochastic Maxwell's equation
\begin{align*}
\EE\left[\Phi(t_k)\right]=\EE\left[\Phi(0)\right]+Kt_k\quad\text{for all discrete times}\quad t_k,
\end{align*}
where we denote $\displaystyle\Phi(t_k):=\int_{\OO}\left(\|\E_k\|^2+\|\H_k\|^2\right)\,\text d{\bf x}$ the numerical energy, 
recall that $t_k=k\Delta t$ for $k= 1,2,\ldots$ and $K=3\left(\lambda_1^2+\lambda_2^2\right)\Tr(Q)$ as in the above result.
\end{prop}
{\noindent\bf Proof.} We first observe that $\Phi(t_k)$ stands for the norm $\|\UU_k\|_V^2$ 
which we now compute
\begin{align*}
\|\UU_k\|^2_V&=\|\S(\Delta t)\UU_{k-1}\|_V^2+2\<\S(\Delta t)\UU_{k-1},\S(\Delta t)\GG\Delta W_{k-1}\>_V\\
&\quad+\|\S(\Delta t)\GG\Delta W_{k-1}\|_V^2\\
&=\|\UU_{k-1}\|_V^2+2\<\S(\Delta t)\UU_{k-1},\S(\Delta t)\GG\Delta W_{k-1}\>_V+\|\GG\Delta W_{k-1}\|_V^2,
\end{align*}
which leads to
\begin{align*}
\EE\left[\|U_k\|^2_V\right]=
\EE\left[\|U_{k-1}\|^2_V\right]+\EE\left[\|\GG\Delta W_{k-1}\|_V^2\right].
\end{align*}
Moreover, using the definition of the $\|\cdot\|_V$ norm and It\^o's isometry, 
one obtains 
\begin{align*}
\EE\left[\|\GG\Delta W_{k-1}\|_V^2\right]
=&3\left(\lambda_1^2+\lambda_2^2\right) \int_\OO 
\EE\left[\left\|\int_{t_{k-1}}^{t_k}\,\text dW(s)\right\|^2\right]\,\text d{\bf x}\\
=&3\left(\lambda_1^2+\lambda_2^2\right)\Delta t \int_\OO
\left(\sum\limits_{n\in \N_+}\eta_ne_n(x)^2\right)
\,\text d{\bf x}\\
=&3\left(\lambda_1^2+\lambda_2^2\right)\Tr(Q)\Delta t=K\Delta t.
\end{align*}
A recursion concludes the proof.\hfill$\square$

The above proposition thus shows that the exact trace formula for the energy also 
holds for the numerical solution given by the exponential integrator \eqref{sm1}.
The following proposition shows that the exponential integrator \eqref{sm1} also 
preserves the discrete version of the averaged divergence exactly.

\begin{prop}\label{div}
The numerical approximation to the linear stochastic Maxwell's equation \eqref{mod;Max1} 
given by the exponential integrator \eqref{sm1}
exactly preserves the following discrete averaged divergence
\begin{align*}
\EE\left[\ddiv (\E_k)\right]=\EE\left[\ddiv (\E_{k-1})\right],\quad 
\EE\left[\ddiv (\H_k)\right]=\EE\left[\ddiv (\H_{k-1})\right]
\end{align*}
for all $k\in\N_+$. 
\end{prop}
{\noindent\bf Proof.}
Let us denote $(\ddiv, \ddiv) (\E^T, \H^T)^T:= (\ddiv \E^T, \ddiv \H^T)^T$.
Taking now the divergence and expectation of both components of 
the numerical solution leads to 
\begin{align}\label{exp-div}
\EE\left[(\ddiv, \ddiv)\UU_{k} \right]=\EE\left[(\ddiv, \ddiv)(\S(\Delta t)\UU_{k-1})\right].
\end{align}
We next notice that $\S(\Delta t)\UU_{k-1}$ is the solution 
of the deterministic Maxwell's equation at time $t=\Delta t$,
\begin{align*}
&\text d\E-\nabla\times\H\,\text dt=0,\nonumber\\
&\text d\H+\nabla\times\E\,\text dt=0, \quad (\E^T,\H^T)^T(0)=\UU_{k-1}.
\end{align*}
Using the property $\ddiv(\nabla\times\cdot)=0$ and a similar argument as 
in \cite[Theorem 2.2]{CHZ16}, we obtain 
\begin{align}\label{lin-div}
(\ddiv, \ddiv)(\S(\Delta t)\UU_{k-1})=(\ddiv, \ddiv )(\UU_{k-1}).
\end{align}
Finally, combing \eqref{exp-div} and \eqref{lin-div} yields the desired result.\hfill$\square$

Regarding the symplectic structure of the numerical solutions, we obtain the following result.  
\begin{prop}\label{symp}
The exponential integrator \eqref{sm1} has the discrete stochastic symplectic conservation law
\begin{align*}
\overline\omega_{1}=
\displaystyle\int_{\OO}\text d\E_{1}\wedge \text d\H_{1}\,\text d{\bf x}=
\displaystyle\int_{\OO}\text d\E_{0}\wedge \text d\H_{0}\,\text d{\bf x}=
\overline\omega_0\quad\mathbb{P}\text{-a.s.}
\end{align*}
\end{prop}
{\noindent\bf Proof.}
Taking the differential of the numerical solution \eqref{sm1} gives 
$\text d \UU_{k+1}=\text d\bigl(\S(\Delta t)\UU_k\bigr)$. Thus, showing symplecticity 
of the exponential integrator is equivalent to showing the symplecticity of the flow of 
the deterministic linear Maxwell's equation with initial value $\UU_k$. 
This is a well know fact.
\hfill$\square$

\section{Numerical experiments}\label{sect-NE}
This section presents various numerical experiments in order to illustrate the main properties of the stochastic 
exponential integrator \eqref{sch;exp}, denoted by \textsc{SEXP} below. 
We will compare this numerical scheme with the following classical ones: 
\begin{itemize}
\item The Euler--Maruyama scheme (denoted by \textsc{EM} below)
\begin{equation}
\label{EM}\tag{EM}
\UU_{k+1}=\UU_k+A\UU_k\Delta t+\FF(\UU_k)\Delta t+\GG(\UU_k)\Delta W_k.
\end{equation}
\item The semi-implicit Euler--Maruyama scheme (denoted by \textsc{SEM} below)
\begin{equation}
\label{SEM}\tag{SEM}
\UU_{k+1}=\UU_k+A\UU_{k+1}\Delta t+\FF(\UU_k)\Delta t+\GG(\UU_k)\Delta W_k.
\end{equation}
\end{itemize}
Below, we consider the stochastic Maxwell's equation \eqref{mod;CMAX} with TM polarization on the domain $[0,1]\times[0,1]$. 
In this setting, the electric and magnetic fields are $\mathbf E=(0,0,E_3)$, resp. $\mathbf H=(H_1,H_2,0)$. 
The spatial discretisation is done by the stagged uniform grid from \cite{v11} with mesh sizes $\Delta x=\Delta y=2^{-4}$. 
Unless stated otherwise, the initial condition reads 
\begin{align*}
E_3(x,y,0)&=0.1\exp(-50((x-0.5)^2+(y-0.5)^2))\\
H_1(x,y,0)&=\text{rand}_y\\
H_2(x,y,0)&=\text{rand}_x,
\end{align*}
where $\text{rand}_x$, resp. $\text{rand}_y$, are random initial values in one direction whereas the other direction is kept constant. 
This is done in order to have zero divergence. The eigenvalues of the linear operator $Q$ are given by $3/(j^3+k^3)$ for $j,k=1,2,\ldots$.

\subsection{Strong convergence}
\begin{figure}[h]
\centering
\includegraphics*[height=4.4cm,keepaspectratio]{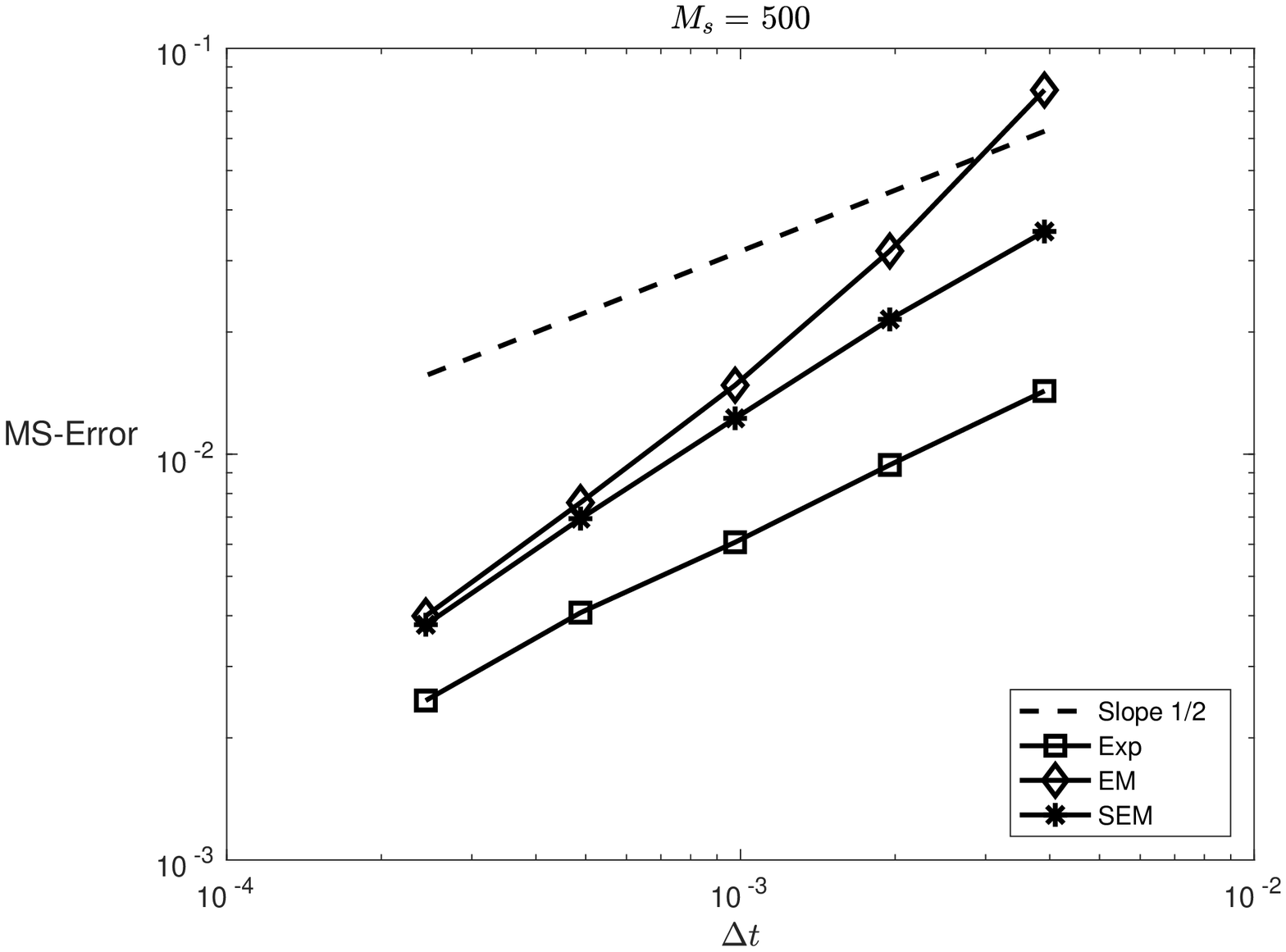}
\includegraphics*[height=4.4cm,keepaspectratio]{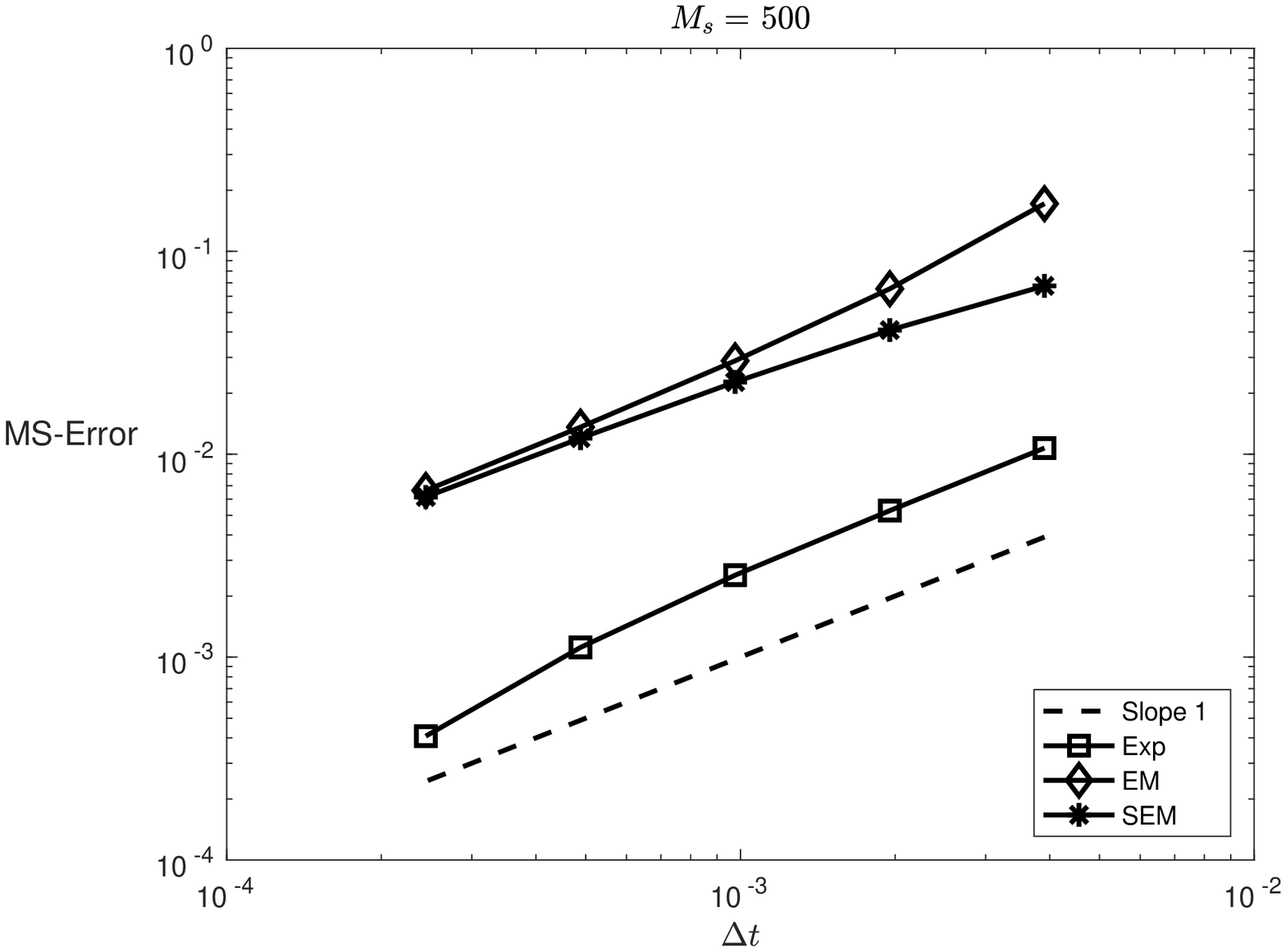}
\caption{Strong rates of convergence for the stochastic Maxwell's equation with 
$\FF(\UU)=\UU+\cos(\UU)$ 
and $\GG(\UU)=\sin(\UU)$ (left) 
and $\FF(\UU)=\UU$ and $\GG(\UU)=\mathbbm{1}^T$ (right).}
\label{fig:strong}
\end{figure}

We first illustrate the strong rates of convergence of the exponential integrator \eqref{sch;exp} 
stated in Theorems~\ref{thm-strong-add}~and~\ref{thm-strong}. To do this, we compute the errors 
$\E\left[\|{\UU^N-\UU_{\text{ref}}(T)}\|^2_{V}\right]$ at the final time $T=0.5$ for time steps ranging from 
$\Delta t=2^{-8}$ to $\Delta t_{\text{ref}}=2^{-13}$ 
and report these errors in Figure~\ref{fig:strong}. 
The reference solution is computed using the exponential integrator 
and the expected values are approximated by computing averages over $M_s=500$ samples. 
We observed that using a larger number of samples ($M_s=750$) does not significantly improve 
the behavior of the convergence plots. The theoretical rates of convergence of the exponential integrator stated in the above theorems 
are indeed observed in these plots. 

\subsection{Averaged energy and divergence}
\begin{figure}[h]
\centering
\includegraphics*[height=2.9cm,keepaspectratio]{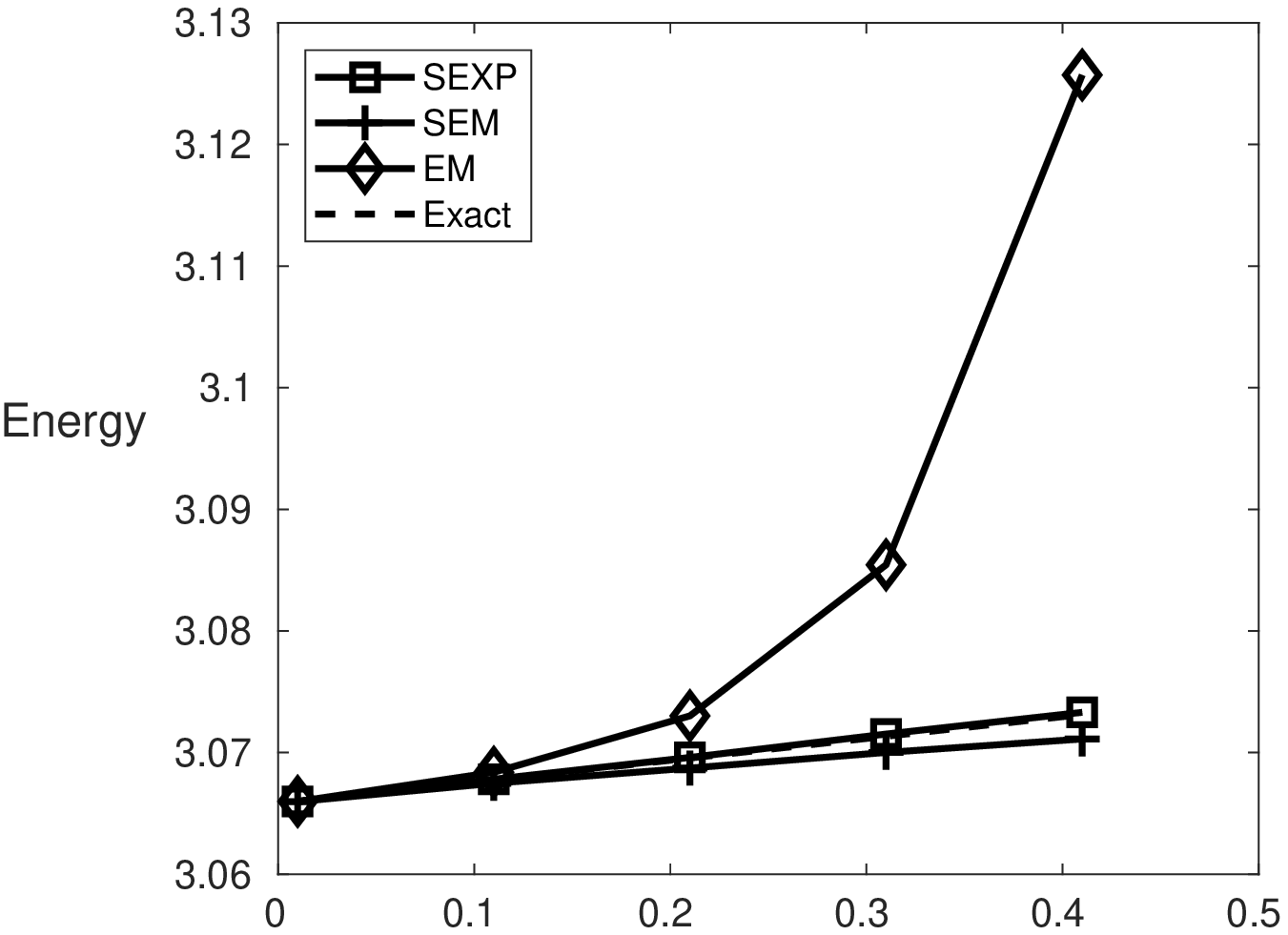}
\includegraphics*[height=2.9cm,keepaspectratio]{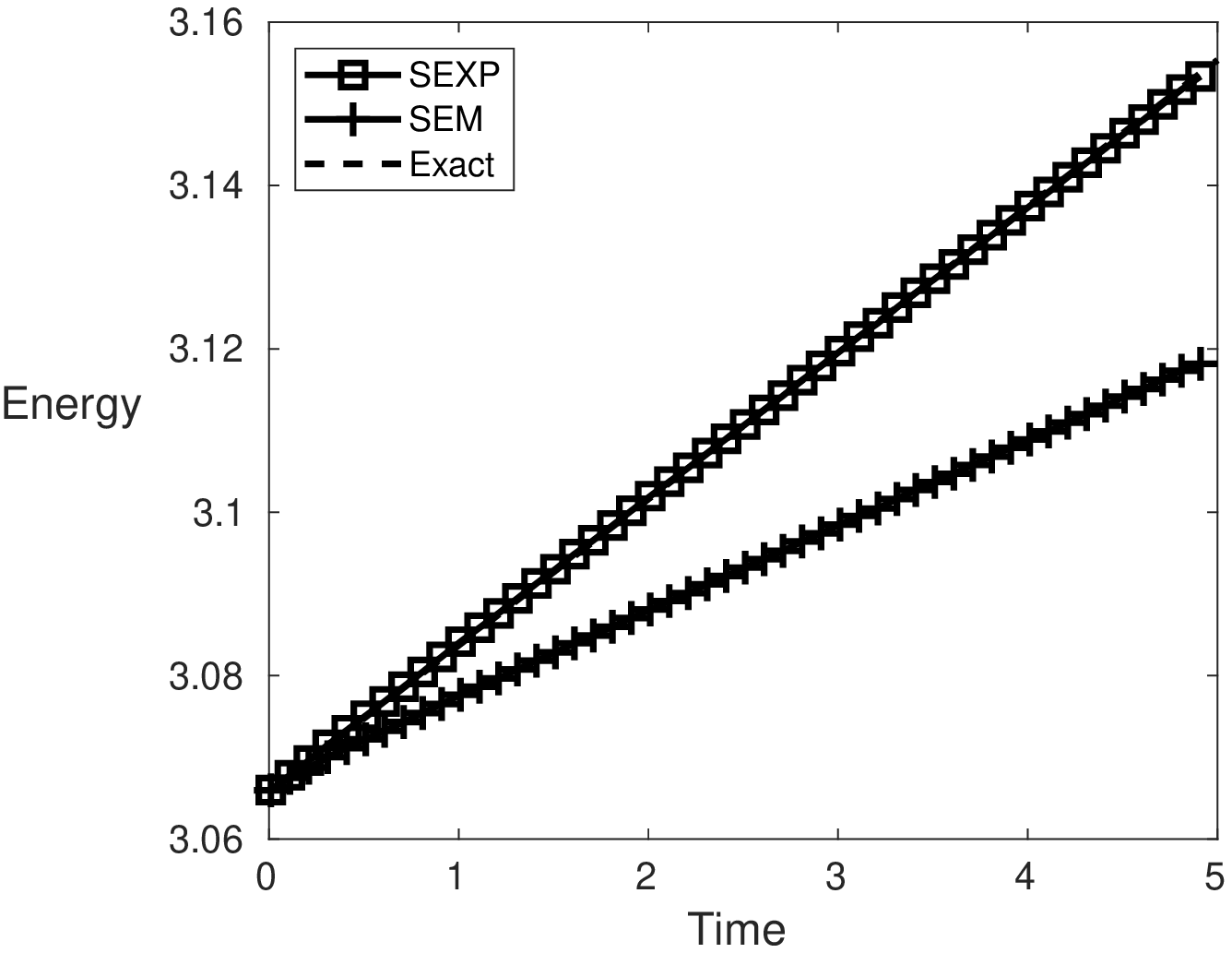}
\includegraphics*[height=2.9cm,keepaspectratio]{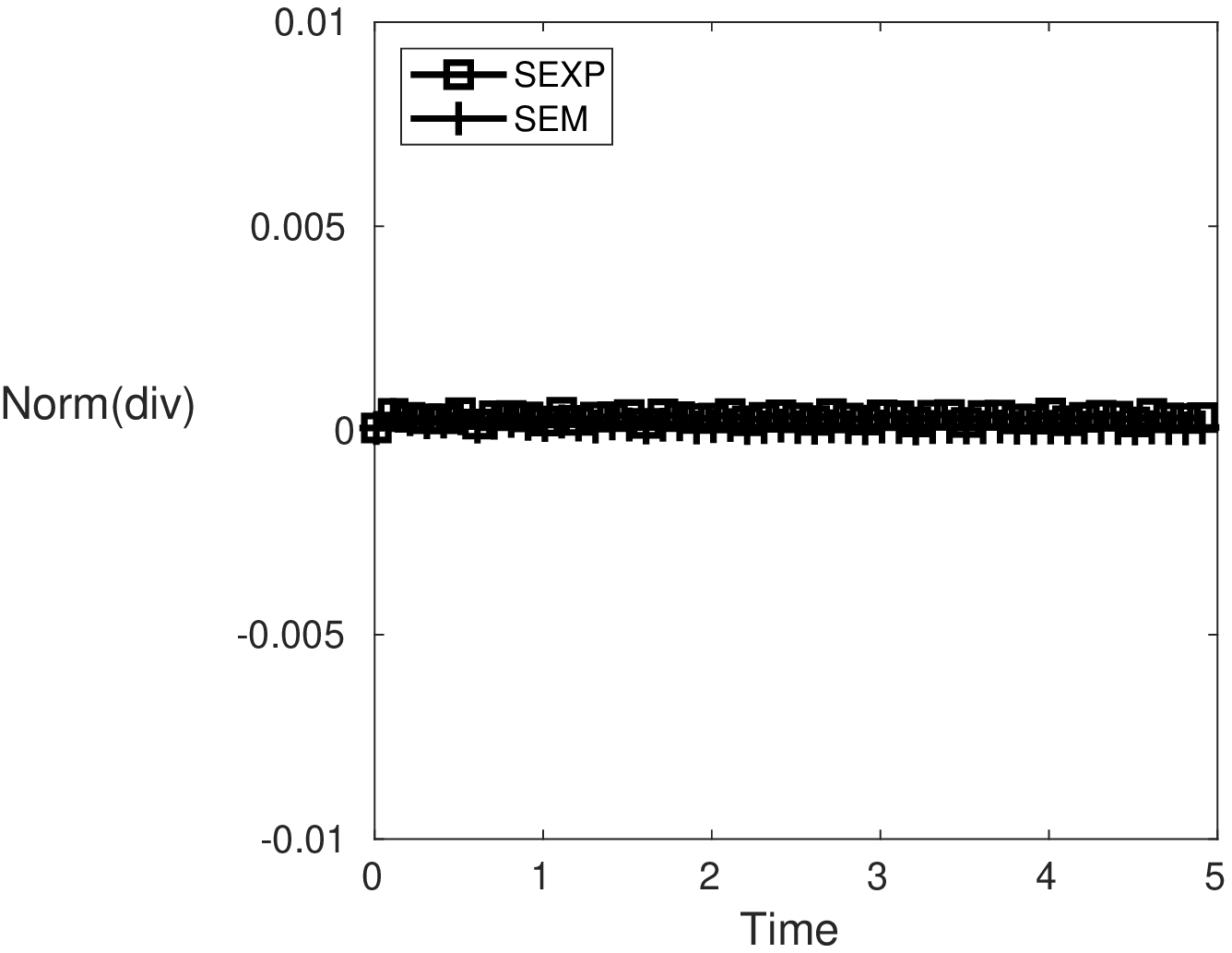}
\caption{Averaged energy on a short time (left) and on a longer time (middle), averaged divergence (right).}
\label{fig:trace}
\end{figure}

We now illustrate the geometric properties of the exponential integrator stated in Section~\ref{sect-LSM}. 
We consider the problem \eqref{mod;Max1} with $\lambda_1=\lambda_2=0.5$, the time interval $[0,5]$, 
a step size $\Delta t=0.01$ and $M_s=25000$ samples to approximate the expectations. 
The numerical averaged energies and divergences are displayed in Figure~\ref{fig:trace}. 
The trace formula for the energy of the stochastic exponential integrator, 
as stated in Proposition~\ref{tm;AvEn}, is observed in this figure (left and middle plots). 
This is in contrast with the wrong behavior of the SEM scheme and the EM scheme, 
where explosion in the energy is observed for the EM scheme (left plot). In this figure (right plot), one can 
also observe the preservation of the averaged divergence of the magnetic field along 
the numerical solution given by the exponential integrator. This confirms the result 
of Proposition~\ref{div}.

\bibliographystyle{amsplain}
\bibliography{bib}

\end{document}